\documentclass[11pt]{amsart}

\theoremstyle{plain}
\newtheorem{theorem}                 {Theorem}      [section]
\newtheorem{proposition}  [theorem]  {Proposition}
\newtheorem{corollary}    [theorem]  {Corollary}

\newtheorem{definition}   [theorem]  {Definition}

\theoremstyle{definition}

\newtheorem{remark}       [theorem]  {Remark}

\numberwithin{equation}{section}

\def \1{\mbox{${\mathbf 1}$}}
\def \r{\mbox{${\mathbb R}$}}
\def \s{\mbox{${\mathbb S}$}}

\def \e{\mbox{${\mathbb E}$}}
\def \h{\mbox{${\mathbb{H}}$}}

\def \tph{\mbox{${\scriptstyle \phi}$}}

\def \f{\mbox{$\varphi$}}

\def \rm{\mbox{${\scriptstyle \frac{1}{\sqrt 2}}$}}
\def \rH{\mbox{${\scriptstyle  \frac{1}{\sqrt{1+|H|^2}}}$}}

\DeclareMathOperator{\trace}{trace}
\DeclareMathOperator{\grad}{grad}

\DeclareMathOperator{\id}{Id}

\DeclareMathOperator{\riem}{Riem}
\DeclareMathOperator{\ricci}{Ricci}

\begin{document}
\title[]
{Classification results for biharmonic submanifolds in spheres}
\author{A. Balmu\c s}
\address{Universit\`a degli Studi di Cagliari\\
Dipartimento di Matematica\\
\newline
Via Ospedale 72\\
09124 Cagliari, ITALIA} \email{balmus@unica.it, montaldo@unica.it}
\author{S. Montaldo}
\author{C. Oniciuc}
\address{Faculty of Mathematics, ``Al.I.~Cuza'' University of Iasi\\
\newline
Bd. Carol I Nr. 11 \\
700506 Iasi, ROMANIA} \email{oniciucc@uaic.ro}

\dedicatory{Dedicated to Professor Vasile Oproiu on his 65th
birthday}


\thanks{The first author was supported by a INdAM doctoral fellowship, ITALY.
The third author was supported by Grant CEEX ET 5871/2006,
ROMANIA}

\begin{abstract}  We classify biharmonic submanifolds with
certain geometric properties in Euclidean spheres. For codimension
$1$, we determine the biharmonic hypersurfaces with at most two
distinct principal curvatures and the conformally flat biharmonic
hypersurfaces. We obtain some rigidity results for
pseudo-umbilical biharmonic submanifolds of codimension $2$ and
for biharmonic surfaces with parallel mean curvature vector field.
We also study the type, in the sense of B-Y.~Chen, of compact
proper biharmonic submanifolds with constant mean curvature in
spheres.
\end{abstract} \maketitle

\section{Introduction}

The study of biharmonic maps between Riemannian manifolds, as a
generalization of harmonic maps, was suggested by J.~Eells and
J.H.~Sampson in \cite{EelSam64}. They define the {\it energy} of a
smooth map $\phi:(M,g)\to(N,h)$ between two Riemannian manifolds,
by $E(\phi)=\frac{1}{2}\int_{M}\, |d\phi|^2\,v_g$, and say that
$\phi$ is {\it harmonic} if it is a critical point of the energy.
The Euler-Lagrange equation associated to $E$ is given by the
vanishing of the {\it tension field} $\tau(\phi)=\trace\nabla
d\phi$.

By integrating the square of the norm of the tension field one can
consider the {\it bienergy} of a smooth map $\phi$,
$E_2(\phi)=\frac{1}{2}\int_{M}\, |\tau(\phi)|^2\,v_g$, and define
its critical points {\it biharmonic} maps (see \cite{BIB}). The
first variation formula for the bienergy, derived in \cite{GYJ1},
shows that the Euler-Lagrange equation associated to $E_2$ is
given by the vanishing of the {\it bitension field} $
\tau_2(\phi)=-J^{\tph}(\tau(\phi))=-\Delta\tau(\phi) -\trace \
R^N(d\phi,\tau(\phi))d\phi$, where $J^{\tph}$ is formally the
Jacobi operator of $\phi$. The operator $J^{\tph}$ is obviously
linear, thus any harmonic map is biharmonic. We call {\it proper
biharmonic} the non-harmonic biharmonic maps.

During the last decade important progress has been made in the
study of  both the geometry and the analytic properties of
biharmonic maps. In differential geometry, a special attention has
been payed to the study of biharmonic submanifolds, i.e.
submanifolds such that the inclusion map is a biharmonic map.


Moreover, the non-existence theorems for the case of non-positive
sectional curvature codo\-mains, as well as the
\\{\bf Generalized Chen's
Conjecture:} {\it Biharmonic submanifolds of a manifold $N$ with
$\riem^N\leq 0$ are minimal}, \\encouraged the study of proper
biharmonic submanifolds in spheres or other non-negatively curved
spaces \cite{CMO1, CMO2, Fet, INO, MO, OU}.

Although important results and examples were obtained, the
classification of proper biharmonic submanifolds in spheres is
still an open problem.

This paper is fully devoted to the classification of proper
biharmonic submanifolds with certain geometric properties in
spheres.

It is organized as follows. In the preliminary section we remind
some fundamental characterization theorems and results on proper
biharmonic submanifolds of space forms and, in particular, of the
Euclidean sphere. This section also contains some basic
information on finite type Euclidean submanifolds. Although
defined in a different manner, finite type submanifolds are, in a
natural way, solutions of a variational problem. They are critical
points of the volume functional for a certain class of directional
deformations (see \cite{BYC}).

In the third section we study the type of compact proper
biharmonic submanifolds of constant mean curvature in $\s^n$ and
prove that they are $1-$type or $2-$type submanifolds of
$\r^{n+1}$.

The fourth section is devoted to the complete classification of the
proper biharmonic hypersurfaces with at most two distinct principal
curvatures in $\s^{m+1}$. We prove that they are open parts of the
hypersphere $\s^m(\rm)$ or of the Clifford tori $\s^{m_1}(\rm)\times
\s^{m_2}(\rm)$, $m_1+m_2=m$, $m_1\neq m_2$ (see Theorem \ref{th:
classif_hypersurf_2_curv_princ}). A similar result is obtained for
conformally flat biharmonic hypersurfaces in spheres. On the
contrary, for the hyperbolic space $\h^{m+1}$ we prove non-existence
results for such hypersurfaces.

In the fifth section we prove that the pseudo-umbilical biharmonic
submanifolds in spheres have constant mean curvature and we give
an estimate for their scalar curvature. Then we classify proper
biharmonic pseudo-umbilical submanifolds of codimension $2$
(Theorem \ref{th: classif_pseudo_umb_codim2}). We also prove that
the only biharmonic surfaces with parallel mean curvature vector
field in $\s^n$ are the minimal surfaces of $\s^{n-1}(\rm)$ (see
Theorem \ref{th: classif_surf_parallelH}).

Also, based on all the known results on biharmonic submanifolds in
spheres, we suggest two conjectures.


\section{Preliminaries}

\subsection{Biharmonic submanifolds}\quad

Consider $\phi:M\to\e^n(c)$ to be the canonical inclusion of a
submanifold $M$ in a constant sectional curvature $c$ manifold,
$\e^n(c)$. The expressions assumed by the tension and bitension
fields are
$$
\tau(\phi)=mH,\qquad\qquad \tau_2(\phi)=-m(\Delta H-mcH),
$$
where $H$ denotes the mean curvature vector field of $M$ in
$\e^n(c)$.

The attempt of classifying the biharmonic submanifolds in space
forms was initiated in \cite{BYC_ISH} and \cite{CMO1} with the
following characterization results, obtained by splitting the
bitension field in its normal and tangent components.

\begin{theorem}\cite{CMO1}.
The canonical inclusion $\phi:M^m\to\e^n(c)$ of a submanifold $M$
in an $n$-dimensional space form $\e^n(c)$ is biharmonic if and
only if
\begin{equation}\label{caract_bih_spheres}
\left\{
\begin{array}{l}
\ -\Delta^\perp H-\trace B(\cdot,A_H\cdot)+mcH=0,
\\ \mbox{} \\
\ 2\trace A_{\nabla^\perp_{(\cdot)}H}(\cdot)
+\frac{m}{2}\grad(\vert H \vert^2)=0,
\end{array}
\right.
\end{equation}
where $A$ denotes the Weingarten operator, $B$ the second
fundamental form, $H$ the mean curvature vector field,
$\nabla^\perp$ and $\Delta^\perp$ the connection and the Laplacian
in the normal bundle of $M$ in $\e^n(c)$.
\end{theorem}

For hypersurfaces, this result becomes

\begin{proposition}\label{prop: caract_hypersurf_bih}
Let $M$ be a hypersurface of $\e^{m+1}(c)$. Then $M$ is proper
biharmonic if and only if
\begin{equation}\label{caract_bih_hipersurf_spheres}
\left\{
\begin{array}{l}
\Delta^\perp H-(mc-|A|^2)H=0,
\\ \mbox{} \\
\ 2A\big(\grad (|H|)\big)+m|H|\grad(|H|)=0.
\end{array}
\right.
\end{equation}
\end{proposition}

In the case of the hyperbolic space some non-existence results
were given. We recall here
\begin{theorem}\cite{CMO2}.
Any biharmonic pseudo-umbilical submanifold $M^m$, $m\neq 4$, of the
hyperbolic space $\h^{n}$ is minimal.
\end{theorem}

For the sphere, using the canonical inclusion in the Euclidean
space, the next caracterization result was obtained
\begin{theorem}\label{th: bitens comp R}\cite{CMO2}.
If $\phi:(M,g)\to\s^n$ is a Riemannian immersion and
$\f=\mathbf{i}\circ\phi$, where
$\mathbf{i}:\s^n\to\r^{n+1}$ is the canonical inclusion,
then
$$
\tau_2(\phi)=\tau_2(\f)+2m\tau(\f)+\{2m^2-|\tau(\f)|^2\}\f.
$$
\end{theorem}

The first achievement towards the classification problem is
represented by the complete classification of proper biharmonic
submanifolds of the $3$-dimensional unit Euclidean sphere,
obtained in \cite{CMO1}.

\begin{theorem}\cite{CMO1}.\label{th: classif_bih_s3}
\begin{itemize}
\item[a)] An arc length parameterized curve
$\gamma:I\to\s^3$ is proper biharmonic if and only if it is
either the circle of radius $\rm$, or a geodesic of the
Clifford torus $\s^1(\rm)\times\s^1(\rm)\subset\s^3$ with
slope different from $\pm 1$.

\item[b)] A surface $M$ is proper biharmonic in $\s^3$ if
and only if it is locally a piece of
$\s^2(\rm)\subset\s^3$. Furthermore, if $M$ is compact and
orientable, then it is proper biharmonic if and only if
$M=\s^2(\rm)$.

\end{itemize}
\end{theorem}

Then, inspired by the $3-$dimensional case, two methods for
constructing proper biharmonic submanifolds in $\s^n$ were given.

\begin{theorem}\cite{CMO2}.\label{th: rm_minim}
Let $M$ be a minimal submanifold of\, $\s^{n-1}(a)\subset\s^n$.
Then $M$ is proper biharmonic in $\s^{n}$ if and only if $a=\rm$.
\end{theorem}

\begin{remark}\quad

\begin{itemize}
\item[a)] This result proved to be quite useful for the construction of proper
biharmonic submanifolds in spheres. For instance, it
implies the existence of closed orientable embedded proper
biharmonic surfaces of arbitrary genus in $\s^4$ (see
\cite{CMO2}).

\item[b)] All minimal submanifolds of
$\s^{n-1}(\rm)\subset\s^n$ are pseudo-umbilical, have pa\-rallel
mean curvature vector in $\s^n$ and $|H|=1 $.
\end{itemize}
\end{remark}

Non pseudo-umbilical examples were also produced by proving

\begin{theorem}\label{th:hipertor}\cite{CMO2}.
Let $M_1^{m_1}$ and $M_2^{m_2}$ be two minimal submanifolds
of $\s^{n_1}(r_1)$ and $\s^{n_2}(r_2)$, respectively, where
$n_1+n_2=n-1$, $r_1^2+r_2^2=1$. Then $M_1\times M_2$ is
proper biharmonic in $\s^n$ if and only if $r_1=r_2=\rm$
and $m_1\neq m_2$.
\end{theorem}

\begin{remark}\quad

\begin{itemize}
\item[a)] The proper biharmonic submanifolds of $\s^n$ constructed as
above are not pseudo-umbilical, but have parallel mean curvature
vector field, thus constant mean curvature, i.e. constant norm of
the mean curvature vector field, and $|H|\in (0,1)$.
\item[b)] The generalized Clifford torus, $\s^{n_1}(\rm)\times\s^{n_2}(\rm)$,
$n_1+n_2=n-1$, $n_1\neq n_2$, was the first example of proper
biharmonic submanifold in $\s^n$ (see \cite{GYJ1}).

\end{itemize}
\end{remark}

We end this section with a partial classification result for
constant mean curvature biharmonic submanifolds in spheres. The
result was obtained in \cite{O2} and due to its importance for our
paper we shall present it with its proof.

\begin{theorem}\label{classif_bih const mean}\cite{O2}.
Let $M$ be a proper biharmonic submanifold with constant mean
curvature $|H|$ in $\s^n$. Then $|H|\in(0,1]$. Moreover, if
$|H|=1$, then $M$ is a minimal submanifold of a hypersphere
$\s^{n-1}(\rm)\subset\s^n$.
\end{theorem}

\begin{proof}

Let $M$ be a constant mean curvature biharmonic submanifold of
$\s^n$. The first equation of \eqref{caract_bih_spheres} implies
that
$$
\langle \Delta^\perp H,H\rangle=m|H|^2-|A_H|^2,
$$
and by using the Weitzenb\" ock formula,
$$
\frac{1}{2}\Delta|H|^2=\langle \Delta^\perp
H,H\rangle-|\nabla^\perp H|^2,
$$
it follows
\begin{equation}\label{expr_bih}
m|H|^2=|A_H|^2+|\nabla^\perp H|^2.
\end{equation}
Let now $\{X_i\}$ be a local orthonormal basis such that
$A_H(X_i)=\lambda_iX_i$. From
$$
\lambda_i=\langle A_H(X_i),X_i\rangle=\langle
B(X_i,X_i),H\rangle
$$
and
$$
\sum \lambda_i=m|H|^2,\qquad \sum(\lambda_i)^2=|A_H|^2,
$$
using \eqref{expr_bih} we obtain
\begin{equation}\label{aux}
\sum \lambda_i=\sum (\lambda_i)^2+|\nabla^\perp H|^2\geq
\frac{(\sum \lambda_i)^2}{m}+|\nabla^\perp H|^2.
\end{equation}
Thus
$$
m|H|^2\geq m|H|^4+|\nabla^\perp H|^2.
$$
Consequently, if $|H|>1$, the last inequality leads to a
contradiction.

If $|H|=1$, then the last inequality implies $\nabla^\perp H=0$
and $\sum (\lambda_i)^2= \frac{(\sum \lambda_i)^2}{m}=m$, thus we
get $\lambda_1=\ldots=\lambda_m$. Therefore $M$ is a minimal
submanifold of the hypersphere $\s^{n-1}(\rm)$.
\end{proof}

\subsection{Pseudo-umbilical submanifolds in spheres}\qquad

\begin{definition}
A submanifold $M$ of a Riemannian manifold $N$ is said to be
pseudo-umbilical if there exists a function $\lambda\in
C^\infty(M)$, such that $A_H=\lambda\id$, where $A_H$ is the
Weingarten operator associated to the mean curvature vector field
$H$ of $M$ in $N$.
\end{definition}

\begin{remark}
If $M$ is a pseudo-umbilical submanifold of $N$, one can
immediately prove that $\lambda=|H|^2$.
\end{remark}

We also recall here two important geometric properties of
pseudo-umbilical submanifolds in spheres.

\begin{theorem}\cite[p.173]{BYC0}\label{th: byc1}.
Let $M$ be an $m-$dimensional pseudo-umbilical submanifold of an
$n-$dimensional unit Euclidean sphere $\s^n$. Then the scalar
curvature $\tau$ of $M$ satisfies
$$
\tau\leq m(m-1)(1+|H|^2).
$$
The equality holds if and only if $M$ is contained in an
$m-$sphere $\s^m\big(\rH\big)$ of $\s^n$.
\end{theorem}

\begin{theorem}\cite[p.180]{BYC0}\label{th: byc2}.
Let $M^m$ be a pseudo-umbilical submanifold in $\s^{m+2}$. If $M$
has constant mean curvature, then $M$ is either a minimal
submanifold of $\s^{m+2}$ or a minimal hypersurface of a
hypersphere of $\s^{m+2}$.
\end{theorem}

\subsection{Finite type submanifolds in Euclidean spaces}

\begin{definition}
An isometric immersion $\varphi:M\to \r^{n}$ is called of finite
type if $\f$ can be expressed as a finite sum of $\r^{n}-$valued
eigenfunctions of the Laplacian $\Delta$ of $M$. When $M$ is
compact it is called of $k-$type if the spectral decomposition of
$\f$ contains exactly $k$ non-zero terms, excepting the center of
mass.
\end{definition}

The following result constitutes a useful tool in determining
whether a compact submanifold of $\r^n$ is of finite type.

\begin{theorem}\label{Chen_crit}(Minimal Polynomial
Criterion).\cite{BYC1, BYC}. Let $\f:M^m\to\r^{n}$ be an isometric
immersion of a compact Riemannian manifold $M$ into $\r^{n}$ and
denote by $H^0$ the mean curvature vector field of $M$ in
$\r^{n}$. Then
\begin{enumerate}
\item[a)] $M$ is of finite type if and only if there exists a
non-trivial polynomial $Q(t)$ such that $Q(\Delta)H^0=0$.

\item[b)] $M$ is of finite type $k$ if and only if there exists a
unique monic (i.e. with leading coefficient equal to $1$)
polynomial $P(t)$ with exactly $k$ distinct positive roots, such
that $P(\Delta)H^0=0$.

\end{enumerate}
\end{theorem}

%
%

%
%

\section{The type of compact proper biharmonic submanifolds in spheres}

In this section, by applying the preliminary results to the
biharmonic case, we intend to analyze the type of proper
biharmonic submanifolds of $\s^n$, as submanifolds in $\r^{n+1}$.
%


We prove the following
%
%

\begin{theorem}\label{tip_subv} Let $M^m$ be a compact constant
mean curvature, $|H|^2=k$, sub\-manifold in  $\s^n$. Then $M$ is
proper biharmonic if and only if

either
\begin{itemize}
\item[a)]  $|H|^2=1$ and $M$ is a 1-type
submanifold of $\r^{n+1}$ with eigenvalue $\lambda=2m$,
\end{itemize}

or
\begin{itemize}
\item[b)] $|H|^2=k\in (0,1)$ and $M$ is a 2-type
submanifold of $\r^{n+1}$ with the eigenvalues
$\lambda_{1,2}=m(1\pm\sqrt k)$.
\end{itemize}

\end{theorem}

\begin{proof}
We directly apply Theorem \ref{th: bitens comp R}. Denote by
$\phi:M\to\s^n$ the inclusion of $M$ in $\s^n$ and by
$\mathbf{i}:\s^n\to\r^{n+1}$ the canonical inclusion. Let
$\varphi:M\to\r^{n+1}$, $\varphi=\mathbf{i}\circ \phi$, be the
inclusion of $M$ in $\r^{n+1}$. Denote by $H$ the mean curvature
vector field of $M$ in $\s^n$ and by $H^0$ the mean curvature
vector field of $M$ in $\r^{n+1}$.

The tension fields of the immersions $\phi$ and $\varphi$ are
related by
\begin{eqnarray*}
\tau(\varphi)=\tau(\phi)-m\varphi
\end{eqnarray*}
and from here it follows that $H^0=H-\varphi$.

Also, from Theorem \ref{th: bitens comp R}, we get that
$\tau_2(\phi)=0$ if and only if
\begin{equation}\label{caract_bih_HH}
\Delta H^0-2mH^0+m(|H|^2-1)\varphi=0.
\end{equation}

There are two situations to be analyzed.

If $|H|^2=1$, then $\Delta H^0-2mH^0=0$, and Theorem
\ref{Chen_crit} implies that $M$ is a $1$-type submanifold in
$\r^{n+1}$ with eigenvalue $\lambda=2m$.

If $|H|^2=k\in (0,1)$, then equation \eqref{caract_bih_HH} implies
\begin{eqnarray*}
0&=&\Delta\Delta H^0-2m \Delta H^0+m(k-1)\Delta
\varphi\\
&=&\Delta\Delta H^0-2m \Delta H^0-m^2(k-1)H^0.
\end{eqnarray*}

The monic polynomial with positive distinct roots described in
Theorem \ref{Chen_crit}, which provides the type of the
submanifold $M$, is
$$
P(\Delta)=\Delta^2-2m\Delta^1-m^2(k-1)\Delta^0,
$$
so $M$ is a $2$-type submanifold with eigenvalues
$\lambda_{1,2}=m(1\pm\sqrt k)$.

For the converse, let first $M$ be a constant mean curvature
$|H|=1$ submanifold of $\s^n$. Suppose it is of $1-$type with
eigenvalue $\lambda=2m$ in $\r^{n+1}$. This means that $\Delta
\varphi=2m \varphi$, and by applying $\Delta$, it implies $\Delta
H^0-2mH^0=0$. From here we see that $M$ satisfies equation
\eqref{caract_bih_HH}, i.e. it is biharmonic in $\s^n$.

When $|H|^2=k\in(0,1)$ and $M$ is a $2-$type submanifold in
$\r^{n+1}$ with eigenvalues $\lambda_{1,2}=m(1\pm\sqrt k)$ we have
$$
\varphi=x_1+x_2,
$$
where $\Delta x_i=\lambda_i x_i$, $i=1,2$. Applying the Laplacian
we obtain
$$
H^0=-\{x_1+x_2+\sqrt k(x_1-x_2)\}=-\varphi-\sqrt k(x_1-x_2)
$$
and
$$
\Delta
H^0=-m\{(k+1)\varphi+2(-\varphi-H^0)\}=-m\{(k-1)\varphi-2H^0\}.
$$
Finally, using \eqref{caract_bih_HH}, $M$ is biharmonic in $\s^n$.
\end{proof}

\begin{remark}
Note that, using Theorem \ref{classif_bih const mean}, we can
conclude that all proper biharmonic submanifolds of $\s^n$ with
$|H|=1$ are $1-$type submanifolds in $\r^{n+1}$, independently on
whether they are compact or not.

\end{remark}

\section{The classification of biharmonic hypersurfaces
with at most two distinct principal curvatures in spheres}

We recall that if $M$ is a proper biharmonic umbilical hypersurface
in $\s^{m+1}$, then it is an open part of $\s^m(\rm)$ and that there
exist no proper biharmonic umiblical hypersurfaces in $\r^{m+1}$ or
in the hyperbolic space $\h^{m+1}$ .

Similarly to the case of the Euclidean space (see \cite{Dim}), the
study of proper biharmonic hypersurfaces with at most two distinct
principal curvatures constitutes the next natural step for the
classification of proper biharmonic hypersurfaces in space forms.

We underline the fact that there exist examples of hypersurfaces
with at most two distinct principal curvatures and non-constant
mean curvature in any space form. In the following we show that,
by adding the hypothesis of biharmonicity, the mean curvature
proves to be constant.

\begin{theorem}\label{th: curb_med_const_2_curv_princ}
Let $M$ be a hypersurface with at most two distinct principal
curvatures in $\e^{m+1}(c)$. If $M$ is proper biharmonic in
$\e^{m+1}(c)$, then it has constant mean curvature.
\end{theorem}

\begin{proof}
If $M$ is umbilical we immediately get to the conclusion.

For $M$ non-umbilical, suppose that $|H|$ is not constant. This,
together with the hypothesis for $M$ to be proper biharmonic with
at most two distinct principal curvatures in $\e^{m+1}(c)$,
implies the existence of an open subset $U$ of $M$, with
\begin{equation}\label{cond_f}
\left\{
\begin{array}{l}
\grad_p f\neq 0,
\\ \mbox{} \\
f(p)>0,\qquad\qquad\qquad\forall p\in U
\\ \mbox{} \\
k_1(p)\neq k_2(p),
\\ \mbox{} \\
m_1, m_2\quad \textrm{constant},
\end{array}
\right.
\end{equation}
where, denoting by $\eta$ the unit section in the normal bundle,
$f$ is the mean curvature function of $U$ in $\e^{m+1}(c)$, i.e.
$H=\frac{1}{m}(\trace A)\eta=f\eta$, and $k_1$, $k_2$ are the
principal curvature functions w.r.t. $\eta$, with multiplicities
$m_1, m_2$.

Under these hypotheses, we shall prove that $f$ is constant on
$U$, contradicting the condition $\grad_p f\neq 0, \forall p\in
U$.

Since $M$ is proper biharmonic in $\e^{m+1}(c)$, from
\eqref{caract_bih_hipersurf_spheres} we have
\begin{equation}\label{caract_bih_spheres_hypersurf}
\left\{
\begin{array}{l}
\ \Delta f=(mc-|A|^2)f,
\\ \mbox{} \\
A(\grad f)=-\frac{m}{2}f\grad f.
\end{array}
\right.
\end{equation}

Consider now $X_1=\displaystyle{\frac{\grad f}{|\grad f|}}$ on
$U$. Then $X_1$ is a principal direction with principal curvature
$k_1=-\frac{m}{2}f$. Suppose that there are $m_1$ principal
directions of principal curvature $k_1$ and $m_2$ principal
directions of principal curvature $k_2\neq k_1$ and recall that
$mf=m_1k_1+m_2k_2$.

We shall use the moving frames method and denote by $X_1,
\{X_i\}_{i=2}^{m_1}, \{X_\alpha\}_{\alpha=m_1+1}^m$ the
orthonormal frame field of principal directions and by
$\{\omega^a\}_{a=1}^m$ the dual frame field of $\{X_a\}_{a=1}^m$
on $U$.

Obviously,
$$
X_i(f)=\langle X_i, \grad f\rangle=|\grad f|\langle
X_i,X_1\rangle=0, \qquad i=2,\ldots,m_1
$$
and analogously $X_\alpha(f)=0$, $\alpha=m_1+1,\ldots, m$, thus
$$
\grad f=X_1(f)X_1.
$$

We write
$$
\nabla X_a=\omega_a^bX_b, \qquad \omega_a^b\in C(T^*U).
$$
From the Codazzi equations for $M$ we get
\begin{eqnarray}\label{eq: Codazzi_hyper1}
X_a(k_b)=(k_a-k_b)\omega_a^b(X_b)
\end{eqnarray}
and
\begin{eqnarray}\label{eq: Codazzi_hyper2}
(k_b-k_d)\omega^d_b(X_a)=(k_a-k_d)\omega^d_a(X_b),
\end{eqnarray}
for distinct $a, b, d=1,\ldots, m$.

We shall show, in the first place, that $m_1=1$.
\\Consider in equation \eqref{eq: Codazzi_hyper1}, $a=1, b=i$. This leads to
$X_1(k_1)=0$, thus $|\grad f|=0$ on $U$ and we have a
contradiction. From here it results that $m_1=1$, thus
$$
k_2=\frac{3m}{2(m-1)}f.
$$

Consider now in \eqref{eq: Codazzi_hyper1}, $a=1$ and $b=\alpha$.
We obtain
\begin{equation}\label{eq: omega_fund}
3X_1(f)=-(m+2)f\omega_1^\alpha(X_\alpha).
\end{equation}

For $a=\alpha$, $b=1$, as $0=X_\alpha(k_1)$, equation \eqref{eq:
Codazzi_hyper1} leads to  $\omega_1^\alpha(X_1)=0$ and we can
write
\begin{equation}\label{eq: omega1a(x1)}
\omega_1^a(X_1)=0,\quad\forall a=1,\ldots,m.
\end{equation}

From \eqref{eq: Codazzi_hyper2}, for $a=1$, $b=\alpha$ and
$d=\beta$, with $\alpha\neq\beta$,
 we get
\begin{equation}\label{eq: omega1b(xa)}
\omega_1^\beta(X_\alpha)=0,\qquad \forall \alpha\neq\beta.
\end{equation}

We now compute
\begin{eqnarray}
\Delta f&=&-\mathrm{div}(\grad f)=-\langle \nabla_{X_1}\grad f,
X_1\rangle-\sum_{\alpha=2}^m\langle\nabla_{X_\alpha}\grad f,
X_\alpha\rangle\\
&=&-X_1\big(X_1(f)\big)-X_1(f)\sum_{\alpha=2}^m
\omega^\alpha_1(X_\alpha)\nonumber
\end{eqnarray}

By using \eqref{eq: omega_fund} we get that
\begin{equation}\label{eq: laplace_curv_med}
f\Delta
f=-fX_1\big(X_1(f)\big)+\frac{3(m-1)}{m+2}\big(X_1(f)\big)^2.
\end{equation}

As $|A|^2=k_1^2+(m-1)k_2^2=\frac{m^2(m+8)}{4(m-1)}f^2$ and $M$ is
biharmonic,
$$
\Delta f=(mc-|A|^2)f=\Big(mc-\frac{m^2(m+8)}{4(m-1)}f^2\Big)f,
$$
and equation \eqref{eq: laplace_curv_med} becomes
\begin{equation}\label{eq: first_fundam_x1}
fX_1\big(X_1(f)\big)-\frac{3(m-1)}{m+2}\big(X_1(f)\big)^2-\frac{m^2(m+8)}{4(m-1)}f^4+mcf^2=0
\end{equation}

We shall now use the Gauss and the Cartan structural equations in
order to obtain other information on $f$. We have

\begin{eqnarray*}
d\omega_1^\alpha&=&-\sum_{a=1}^m\omega_1^a\wedge
\omega_\alpha^a-(k_1k_2+c)\omega^1\wedge\omega^\alpha,
\end{eqnarray*}
thus, using equations \eqref{eq: omega1a(x1)} and  \eqref{eq:
omega1b(xa)}, we get
\begin{equation}\label{eq: domega1a}
d\omega_1^\alpha(X_1,X_\alpha)=
-k_1k_2-c=\frac{3m^2}{4(m-1)}f^2-c.
\end{equation}

On the other hand from \eqref{eq: omega1a(x1)} and \eqref{eq:
omega1b(xa)} we obtain
$\omega_1^\alpha=\omega_1^\alpha(X_\alpha)\omega^\alpha$, thus
\eqref{eq: omega_fund} implies
\begin{equation}\label{eq: omega1a}
3X_1(f)\omega^\alpha=-(m+2)f\omega_1^\alpha.
\end{equation}
By differentiating \eqref{eq: omega1a} we obtain
\begin{equation}\label{eq: deriv_curb}
3d\big(X_1(f)\big)\wedge\omega^\alpha+3X_1(f)d\omega^\alpha
=-(m+2)(df\wedge\omega_1^\alpha+fd\omega_1^\alpha).
\end{equation}

We use \eqref{eq: domega1a}, substitute
$$
\big(dX_1(f)\wedge\omega^\alpha\big)(X_1,X_\alpha)=X_1\big(X_1(f)\big),
$$
$$
d\omega^\alpha(X_1,X_\alpha)
=\omega_1^\alpha(X_\alpha),
$$
$$
(df\wedge\omega_1^\alpha)(X_1,X_\alpha)
=X_1(f)\omega_1^\alpha(X_\alpha)
$$
in \eqref{eq: deriv_curb} and obtain
\begin{equation}\label{eq: second_fundam_x1}
fX_1\big(X_1(f)\big)-\frac{m+5}{m+2}\big(X_1(f)\big)^2+
\frac{m^2(m+2)}{4(m-1)}f^4-\frac{m+2}{3}cf^2=0.
\end{equation}

Consider now an arbitrary integral curve $\gamma$ of $X_1$ and
denote by $f'$ and $f''$ the first and the second derivatives of
$f$ along this curve. Equations \eqref{eq: first_fundam_x1} and
\eqref{eq: second_fundam_x1} become, respectively,
\begin{equation}\label{eq: first_fundam}
ff''-\frac{3(m-1)}{m+2}(f')^2-\frac{m^2(m+8)}{4(m-1)}f^4+mcf^2=0
\end{equation}
and
\begin{equation}\label{eq: second_fundam}
ff''-\frac{m+5}{m+2}(f')^2+
\frac{m^2(m+2)}{4(m-1)}f^4-\frac{m+2}{3}cf^2=0,
\end{equation}
along $\gamma$.

Multiplying by $(m+5)$ equation \eqref{eq: first_fundam} and by
$-3(m-1)$ equation \eqref{eq: second_fundam} and summing up, we
get
\begin{equation}\label{eq: prel_f}
(4-m)ff''=\frac{m^2(m^2+4m+9)}{2(m-1)}f^4-(m^2+3m-1)cf^2.
\end{equation}
For $m=4$, equation \eqref{eq: prel_f} implies $f=$constant and
thus drives to the contradiction. For $m\neq 4$, we multiply
equation \eqref{eq: prel_f} by $f'/f$, integrate the result and
obtain
\begin{equation}\label{eq: derivata_1}
(f')^2=\frac{m^2(m^2+4m+9)}{8(4-m)(m-1)}f^4
-\frac{(m^2+3m-1)}{2(4-m)}cf^2+C.
\end{equation}

On the other hand, multiplying by $-1$ equation \eqref{eq:
first_fundam} and adding it to equation \eqref{eq: second_fundam}
leads to
\begin{equation}\label{eq: derivata_2}
(f')^2=\frac{m^2(m+5)(m+2)}{4(4-m)(m-1)}f^4
-\frac{(2m+1)(m+2)}{3(4-m)}cf^2
\end{equation}

From \eqref{eq: derivata_1} and \eqref{eq: derivata_2} we conclude
that $f$ is the solution of a polynomial equation, thus $f$ is
constant along $\gamma$. As $\gamma$ was an arbitrary integral
curve for $X_1$ we have $X_1(f)=0$ on U, thus we arrive to a
contradiction.
\end{proof}

To strengthen the Generalized Chen's Conjecture, as an immediate
consequence of Theorem \ref{th: curb_med_const_2_curv_princ}, we
have the following non-existence result.

\begin{theorem}
There exist no proper biharmonic hypersurface with at most two
distinct principal curvatures in $\h^{m+1}$.
\end{theorem}
\begin{proof}
Suppose that $M$ is a proper biharmonic hypersurface with at most
two distinct principal curvatures in $\h^{m+1}$. From Theorem
\ref{th: curb_med_const_2_curv_princ}, the mean curvature of $M$
is constant, and applying Proposition \ref{prop:
caract_hypersurf_bih} we obtain $|A|^2=-m$ and we conclude.
\end{proof}

The case of the sphere is essentially different. Theorem \ref{th:
curb_med_const_2_curv_princ} proves to be the main ingredient for
the following complete classification of proper biharmonic
hypersurfaces with at most two distinct principal curvatures.

\begin{theorem}\label{th: classif_hypersurf_2_curv_princ}
Let $M^m$ be a proper biharmonic hypersurface with at most two
distinct principal curvatures in $\s^{m+1}$. \\Then $M$ is an open
part of $\s^{m}(\rm)$ or of $\s^{m_1}(\rm)\times \s^{m_2}(\rm)$,
$m_1+m_2=m$, $m_1\neq m_2$.
\end{theorem}

\begin{proof}

By Theorem  \ref{th: curb_med_const_2_curv_princ}, the mean
curvature of $M$ in $\s^{m+1}$ is constant and, by using Proposition
\ref{prop: caract_hypersurf_bih}, we obtain $|A|^2=m$. These imply
that $M$ has constant principal curvatures.

For $|H|^2=1$ we conclude that $M$ is an open part of $\s^m(\rm)$.

For $|H|^2\in (0,1)$ we deduce that $M$ has two distinct constant
principal curvatures. Proposition 2.5 in \cite{RY} implies that $M$
is an open part of the product of two spheres $\s^{m_1}(a)\times
\s^{m_2}(b)$, such that $a^2+b^2=1$, $m_1+m_2=m$. Since $M$ is
biharmonic in $\s^n$, from Theorem \ref{th:hipertor} it follows that
$a=b=\rm$ and $m_1\neq m_2$.
\end{proof}

\begin{remark}
Note that, for $m=2$ we recover the result in Theorem \ref{th:
classif_bih_s3} $b)$.
\end{remark}

We recall that a Riemannian manifold is called conformally flat if
for every point it admits an open neighborhood conformally
diffeomorphic to an open set of an Euclidean space. Also, a
hypersurface $M^m\subset N^{m+1}$ which admits a principal
curvature of multiplicity at least $m-1$ is called
quasi-umbilical.


\begin{theorem}
Let $M^m$, $m\geq 3$, be a  proper biharmonic hypersurface in
$\s^{m+1}$. The following statements are equivalent
\begin{itemize}
\item[a)] $M$ is quasi-umbilical,
\item[b)] $M$ is conformally flat,
\item[c)] $M$ is an open part of $\s^m(\rm)$ or of $\s^{1}(\rm)\times \s^{m-1}(\rm)$.
\end{itemize}
\end{theorem}

\begin{proof}
By Theorem \ref{th: classif_hypersurf_2_curv_princ} we get that
$a)$ is equivalent to $c)$. Also, note that $c)$ obviously implies
$b)$.

In order to prove that $b)$ implies $a)$, remind that, for $m\geq
4$, by a well-known result (see \cite{BYC0}), any conformally flat
hypersurface of a space form is quasi-umbilical and we conclude.

For $m=3$, as $M$ is conformally flat, it results that the
$(0,2)-$tensor field
$\displaystyle{L=-\ricci+\frac{s}{4}\langle\,\,,\,\rangle}$, where
$s$ is the scalar curvature of $M$, is a Codazzi tensor field,
i.e.
\begin{equation}\label{eq: 1}
(\nabla_X L)(Y,Z)=(\nabla_Y L)(X,Z),\quad\forall X, Y, Z\in C(TM).
\end{equation}

Using the notations from the proof of Theorem \ref{th:
curb_med_const_2_curv_princ}, the Gauss equation implies
$$
\ricci(X,Y)=2\langle X,Y\rangle+3f\langle A(X),Y\rangle-\langle
A(X),A(Y)\rangle
$$
and
\begin{equation}\label{eq: 2}
s=6+9f^2-|A|^2.
\end{equation}

We use the same techniques as in the proof of Theorem \ref{th:
curb_med_const_2_curv_princ}. Suppose the existence of an open
subset $U$ of $M$ with $3$ distinct principal curvatures.

If $f$ is constant on $U$, using the above expressions, we
conclude that $U$ is flat and the product of any of its two
principal curvatures is $-1$, thus we get to a contradiction.

Assume that $f$ is not constant on $U$. We can suppose that
$\grad_pf\neq 0$, $\forall p\in U$. Consider
$\displaystyle{X_1=\frac{\grad f}{|\grad f|}}$. As $M$ is proper
biharmonic, $X_1$ gives a principal direction with principal
curvature $\displaystyle{k_1=-\frac{3}{2}f}$. From $k_1+k_2+k_3=3f$,
we can write $\displaystyle{k_2=\frac{9}{4}f+\varepsilon}$ and
$\displaystyle{k_3=\frac{9}{4}f-\varepsilon}$, $\varepsilon\in
C^\infty(U)$. Using the Codazzi and Gauss equations and equations
\eqref{eq: 1} and \eqref{eq: 2} we show that $f=a\varepsilon^5$,
$a\in \r$, and combining all these relations we obtain $\varepsilon$
as the solution of a polynomial equation and we get to a
contradiction.

Finally, it results that $M$ has at most two distinct principal
curvatures and we conclude.
\end{proof}

For what concerns proper biharmonic hypersurfaces with constant
mean curvature in spheres we also have the following geometric
property
\begin{proposition}\label{th: curb_scal_hyp}
Let $M$ be a proper biharmonic hypersurface with constant mean
curvature $|H|^2=k$ in $\s^{m+1}$. Then $M$ has constant scalar
curvature,
$$
s=m^2(1+k)-2m.
$$
\end{proposition}
\begin{proof}
Since $M$ is proper biharmonic of constant mean curvature, the
squared norm of its second fundamental form is $|A|^2=m$. By
applying the Gauss equation, we conclude.
\end{proof}

In view of the above results we propose the following

\subsubsection*{\bf{Conjecture:}} {\it The only proper biharmonic hypersurfaces
in $\s^{m+1}$ are the open parts of hyperspheres $\s^{m}(\rm)$ and
of generalized Clifford tori $\s^{m_1}(\rm)\times \s^{m_2}(\rm)$,
$m_1+m_2=m$, $m_1\neq m_2$.}

\section{Codimension $2$ biharmonic
pseudo-umbilical submanifolds in spheres}

We shall first prove a general result concerning the mean
curvature of biharmonic pseudo-umbilical submanifolds in spheres

\begin{theorem}\label{th: pseudo_bih}
Let $M$ be a pseudo-umbilical submanifold of $\s^{n}$,
$m\neq 4$. If $M$ is biharmonic, then it has constant mean
curvature.
\end{theorem}

\begin{proof}
Consider $x\in M$ and let $\{X_i\}_{i=\overline{1,m}}$ be a local
orthonormal frame field geodesic in $x$. As $M$ is biharmonic,
from \eqref{caract_bih_spheres} we get
\begin{equation}\label{eq: first_pseudoumb}
\trace A_{\nabla^\perp_{({\cdot})}H}
(\cdot)=-\frac{m}{4}\grad(|H|^2).
\end{equation}
On the other hand, in $x$, by standard computations, we get

\begin{eqnarray*}
\trace A_{\nabla^\perp_{({\cdot})}H}(\cdot)
&=&\sum_{i,j}\big\{ X_i\langle
\nabla^{\s^n}_{X_j}X_i,H\rangle- \langle\nabla^{\s^n}_{X_i}
\nabla^{\s^n}_{X_j}X_i,H\rangle\big\}X_j,
\end{eqnarray*}

\begin{eqnarray*}
\sum_{i,j} X_i\langle \nabla^{\s^n}_{X_j}X_i,H\rangle
X_j
&=&\sum_i\nabla_{X_i}A_H(X_i),
\end{eqnarray*}

\begin{eqnarray*}
\sum_{i,j}\langle\nabla^{\s^n}_{X_i}
\nabla^{\s^n}_{X_j}X_i,H\rangle X_j 
&=&\frac{m}{2}\grad(|H|^2).
\end{eqnarray*}

Now, as $M$ is pseudo-umbilical,
$\displaystyle{\sum_i\nabla_{X_i}A_H(X_i)=\grad(|H|^2)}$, thus
\begin{equation}\label{eq: second_pseudoumb}
\trace A_{\nabla^\perp_{(\cdot)}H}(\cdot)=\frac{2-m}{2}\grad
(|H|^2).
\end{equation}

By putting together expressions \eqref{eq: first_pseudoumb} and
\eqref{eq: second_pseudoumb}, we conclude.
\end{proof}

The first consequence of this result is an estimate for the scalar
curvature of biharmonic pseudo-umbilical submanifolds in spheres.

\begin{proposition}\label{th: curb_scal_part}
Let $M^m$ be a biharmonic pseudo-umbilical submanifold of $\s^n$,
$m\neq 4$. Then its scalar curvature $s$ satisfies
$$
s\leq 2m(m-1).
$$
The equality holds if and only if $M$ is open in
$\s^m(\rm)$.
\end{proposition}
\begin{proof}
From Theorem \ref{th: byc1}, the scalar curvature $s$ of a
pseudo-umbilical submanifold $M$ of $\s^n$ satisfies $s\leq
m(m-1)(1+|H|^2)$, and equality holds if and only if $M$ is
contained in an $m$-sphere of $\s^n$.

By using Theorem \ref{th: pseudo_bih} and Theorem \ref{classif_bih
const mean}, as $M$ is biharmonic, it follows that its constant
mean curvature satisfies $|H|\in(0,1]$, and this completes the
proof.
\end{proof}

For what concerns biharmonic pseudo-umbilical submanifolds of
codimension two we obtain the following rigidity result.

\begin{theorem}\label{th: classif_pseudo_umb_codim2}
Let $M^m$ be a pseudo-umbilical submanifold of $\s^{m+2}$, $m\neq
4$. Then $M$ is proper biharmonic if and only if it is minimal in
$\s^{m+1}(\rm)$.
\end{theorem}
\begin{proof}
From the hypotheses, using Theorem \ref{th: pseudo_bih}, we deduce
that $M$ has constant mean curvature. Now, by using Theorem
\ref{th: byc2}, it follows that any such submanifold is either
minimal in $\s^{m+2}$ or minimal in a hypersphere of $\s^{m+2}$.
But $M$ is proper biharmonic in $\s^{m+2}$ and, from Theorem
\ref{th: rm_minim}, we conclude.
\end{proof}

Replace now the condition on $M$ to be pseudo-umbilical with that
of being a hypersurface of a hypersphere in $\s^{m+2}$.

\begin{theorem}\label{th: codim 2, raza a}
Let $M^m$ be a hypersurface of $\s^{m+1}(a)\subset\s^{m+2}$, $a\in
(0,1)$. Assume that $M$ is not minimal in $\s^{m+1}(a)$. Then it is
biharmonic in $\s^{m+2}$ if and only if $a>\rm$ and $M$ is open in
$\s^m(\rm)\subset\s^{m+1}(a)$.
\end{theorem}
\begin{proof}
Note that the converse follows immediately from Theorem \ref{th:
rm_minim}.

In order to prove the other implication, denote by $\mathbf{j}$
and $\mathbf{i}$ the inclusion maps of $M$ in $\s^{m+1}(a)$ and of
$\s^{m+1}(a)$ in $\s^{m+2}$, respectively.

We consider
$$
\s^{m+1}(a)=\Big\{(x^1,\ldots, x^{m+2},\sqrt{1-a^2})\in\r^{m+3}:
\sum_{i=1}^{m+2}(x^i)^2=a^2\Big\}\subset\s^{m+2}.
$$
Then
$$
C\big(T\s^{m+1}(a)\big)=\Big\{(X^1,\ldots, X^{m+2},0)\in
C(T\r^{m+3}): \sum_{i=1}^{m+2}x^iX^i=0\Big\},
$$
while
$$
\eta=\frac{1}{c}\Big(x^1,\ldots,
x^{m+2},-\frac{a^2}{\sqrt{1-a^2}}\Big)
$$
is a unit section in the normal bundle of $\s^{m+1}(a)$ in
$\s^{m+2}$, where $c^2=\frac{a^2}{1-a^2}, c>0$.

By computing the tension and bitension fields of
$\phi=\mathbf{i}\circ \mathbf{j}$, one gets
$$
\tau(\phi)=\tau(\mathbf{j})-\frac{m}{c}\eta,
$$
and
$$
\tau_2(\phi)=\tau_2(\mathbf{j})-\frac{2m}{c^2}\tau(\mathbf{j})+
\frac{1}{c}\big\{|\tau(\mathbf{j})|^2-\frac{m^2}{c^2}(c^2-1)\big\}\eta.
$$

By the hypotheses $M$ is biharmonic in $\s^{m+2}$, thus
$$
|\tau(\mathbf{j})|^2=\frac{m^2}{c^2}(c^2-1)=\frac{m^2}{a^2}(2a^2-1)
$$
and, since $\tau(\mathbf{j})\neq 0$, this implies $a>\rm$.

Also,
$$
|\tau(\phi)|^2=|\tau(\mathbf{j})|^2+\frac{m^2}{c^2}=m^2.
$$
This implies that the mean curvature of $M$ in $\s^{m+2}$ is $1$,
thus, using Theorem \ref{classif_bih const mean}, $M$ has to be a
minimal submanifold of the hypersphere
$\s^{m+1}(\rm)\subset\s^{m+2}$, i.e. it is pseudo-umbilical and
with parallel mean curvature vector field in $\s^{m+2}$.

Since $M\subset \s^{m+1}(a)$ is pseudo-umbilical in $\s^{m+2}$ it
results pseudo-umbilical, and thus totally umbilical, in
$\s^{m+1}(a)$. From here follows that $M$ is an open subset of a
hypersphere $\s^m(r)$ in $\s^{m+1}(a)$. But it is proper
biharmonic in $\s^{m+2}$, thus $r=\rm$ and we conclude.
\end{proof}

\begin{corollary}
Let $M$ be a proper biharmonic hypersurface of a hypersphere
$\s^{m+1}(a)$ in $\s^{m+2}$, $a\in (0,1)$. Then $a\geq\rm$.
Moreover,
\begin{itemize}
\item[a)] if $a=\rm$, then $M$ is minimal in $\s^{m+1}(\rm)$
\item[b)] if $a>\rm$, then $M$ is an open part of $\s^{m}(\rm)$.
\end{itemize}
\end{corollary}

We also use Theorem \ref{th: codim 2, raza a} in order to prove

\begin{theorem}\label{th: classif_surf_parallelH}
Let $M^2$ be a proper biharmonic surface with parallel mean
curvature vector field in $S^n$. Then $M$ is minimal in
$S^{n-1}(\rm)$.
\end{theorem}

\begin{proof}

B-Y.~Chen and S-T.~Yau proved (see \cite[p.106]{BYC0}) that the
only non-minimal surfaces with parallel mean curvature vector
field in $\s^n$ are either minimal surfaces of small hyperspheres
$\s^{n-1}(a)$ of $\s^n$ or surfaces with constant mean curvature
in $3-$spheres of $\s^n$.

If $M$ is a minimal surface of a small hypersphere $\s^{n-1}(a)$,
then it is biharmonic in $\s^n$ if and only if $a=\rm$ (see
Theorem \ref{th: rm_minim}).

If $M$ is a surface in a $3-$sphere $\s^3(a)$, $a\in(0,1]$, of
$\s^n$ then we can consider
$$
M\longrightarrow \s^3(a)\longrightarrow\s^4\longrightarrow\s^n.
$$
Note that $M$ is biharmonic in $\s^n$ if and only if it is
biharmonic in $\s^4$. From Theorem~\ref{th: codim 2, raza a}, for
$a\in(0,1)$, we conclude that either $a=\rm$ and $M$ is minimal in
$\s^3(\rm)$, or $a>\rm$ and $M$ is an open part of $\s^2(\rm)$.
For $a=1$, from Theorem \ref{th: classif_bih_s3}, also follows
that $M$ is an open part of $\s^2(\rm)$.

In all cases $M$ is minimal in $\s^{n-1}(\rm)$.
\end{proof}

\begin{remark}
All the results we have proved so far could suggest that the
codimension $2$ biharmonic submanifolds of $\s^n$ come from minimal
submanifolds of $\s^{n-1}(\rm)$. This is not the case as shown by
the following

\begin{theorem}\cite{S2}.
Let $\phi:M^3\to\s^5$ be a proper biharmonic anti-invariant
immersion. Then the position vector field
$x_0=\mathbf{i}\circ\phi=x_0(u,v,w)$ of $M$ in $\r^6$ is given by

$$
x_0(u,v,w)=\rm e^{iw}(e^{iu},ie^{-iu}\sin\sqrt 2 v,
ie^{-iu}\cos\sqrt 2 v),
$$
where $\mathbf{i}:\s^5\to\r^6$ is the canonical inclusion.
\end{theorem}

Remind that if we consider a Sasakian manifold $(N,
\Phi,\xi,\eta,g)$ and a submanifold $M$ tangent to $\xi$, $M$ is
called {\it anti-invariant} if $\Phi$ maps any tangent vector to
$M$ which is normal to $\xi$ to a vector which is normal to $M$.
Also, a map $\phi:M\to \s^n$ is said to be {\it full} if the image
$\phi(M)$ is contained in no hypersphere of $\s^n$.

Note that $\phi$ is a full proper biharmonic anti-invariant
immersion from a $3-$dimen\-sional torus into $\s^5$. The
immersion $\phi$ has constant mean curvature, is non
pseudo-umbilical and its mean curvature vector is not parallel in
$\s^5$. In addition to these properties, since $|H|=1/3$, from
Theorem \ref{tip_subv} we conclude that $x_0$ is a $2-$type
submanifold of $\r^6$ with eigenvalues $2$ and $4$.

We also note that the product $\s^1(\rm)\times M^m$, where $M$ is
a minimal non-geodesic hypersurface of $\s^{m+1}(\rm)$, is a full
proper biharmonic submanifold of $\s^{m+3}$ of codimension $2$.
\end{remark}

Since all the known examples of proper biharmonic submanifolds in
spheres have constant mean curvature we propose the following

\subsubsection*{\bf{Conjecture:}} {\it Any biharmonic submanifold in
$\s^n$ has constant mean curvature.}

\end{document}